# Innovative Integrated Congestion Management Market: Synergizing Operational Flexibility of Transmission System Operators with Aggregator-Driven Resources

Arash Beiranvand, *Senior Member, IEEE*

***Abstract*—This paper introduces an innovative congestion management market that seeks to integrate and harness the operational flexibility available to Transmission System Operators (TSO) in conjunction with the flexibility offered by aggregators. In the context of operational flexibility, this study investigates Network Topology Optimization (NTO) as a pivotal strategy for TSO to address the transmission network congestion, specifically through the strategic switching of certain buses. To achieve this, an N-1 secured NTO is formulated using a Mixed Integer Linear Programming (MILP) methodology. Concurrently, acknowledging the crucial and growing importance of aggregator-provided flexibility in congestion management, the paper presents a novel, integrated congestion management market. In this framework, aggregators submit bids, and the TSO, taking into account these bids alongside the potential for NTO, clears the market with the goal of minimizing the overall cost of congestion management. The proposed integrated market is constructed using a linear optimization model, ensuring its practical applicability in real-world scenarios. The effectiveness and efficiency of the proposed market design are assessed using the standard IEEE 39-bus system, with results indicating promising potential and applicability.**

## Nomenclature

### Constants

$C^{(m)}_{Gof}$ Generation Adjustment Offer Price (€/MWh) at bus m

$C^{(k)}_{Aof}$ Generation Adjustment 'All-or-Nothing' Offer Price (€/MWh) at Bus $k$

$C^{(w)}_{DGof}$ Generation Adjustment Offer Price for MCP Aggregators (€/MWh) at Bus $w$

$C^{(h)}_{sw}$ The Cost of Switching Operation at Bus $h$

$M$ Sufficiently Large Number

$G^{(i)}_{max}$ Maximum Generation Capacity at Bus $i$ (MW)

$D^{(i)}$ Demand at Bus $i$ (MW)

$P^{(i)}_{Aof}$ Generation Adjustment 'All-or-Nothing' (MW) at Bus $i$

$P^{(i)}_{DGof}$ Generation Adjustment for MCP Aggregators (MW) at Bus $i$

$P^{max}_{l}$ Maximum Loading Capacity of Line $l$ (MW)

$\delta_{max}$ Maximum Allowed Angle

$X_l$ The Reactance of Line $l$

$F_b$ Maximum Number of Allowed Bus Switching Operations

$F_l$ Maximum Number of Allowed Line Switching Operations

### Sets and Indices

$L$ Set of transmission grid's lines

$\hat{L}$ Set of transmission grid's non-radial lines

$B$ Set of grid's buses

$B^{'}_{sw}$ Set of Network Buses Where Switching Is Not Permitted

$L^{'}$ Set of grid's lines except line $k$

$G$ Set of grid's generators

$G^{of}$ Set of Dispatchable Generators Provided by Aggregators With SCP

$G^{Aof}$ Set of 'All-or-Nothing' Flexibilities Provided by Aggregators With SCP

$D^{DG}$ Set of 'All-or-Nothing' Flexibilities Provided by Aggregators With MCP

### Variables and functions

$P^{(m)}_{Gof_1}$ The Offered Generation Adjustments by the SPC Aggregator at Bus $m$, Side 1

$z^{(k)}_{Aof_1}$ Binary Decision Variable for the 'All-or-Nothing' Adjustment Offer at Bus $k$ to be Connected to Side 1, 1 Indicates the Offer Is Accepted

$z^{(w)}_{DGof_1}$ Binary Decision Variable for the MCP Adjustment Offer at Bus $w$ to be Connected to Side 1, 1 Indicates the offer Is Accepted

$P^{(i)}_{l_1}$ Active Power Flow of Line $l$ Connected to Side 1 of Bus $i$ (MW)

$z^{(i)}_{g}$ Binary Decision Variable for the Generator Connected to Bus $i$, 1 Indicates the Generator Is Connected to Side 1

$z^{(i)}_{d}$ Binary Decision Variable for the Demand Connected to Bus $i$, 1 Indicates the Generator Is Connected to Side 1

$z^{(i)}_{of_1}$ Binary Decision Variable for the Offered Generation Adjustment a Bus $i$ to be Connected to Side 1, 1 Indicates the Offer Is Accepted

$z^{(f)}_{l}$ Binary Decision Variable for Line $l$ at Its From-End,

1 Indicates It Is Connected to Side 1

$P_l$ Active Power Flow of Line $l$

$U_l$ Binary Decision Variable for Line $l$ Status, 1 Indicates It Is Out

$\delta_1^{(f)}$ Bus Angle of Side 1 at the From-End of Line $l$

$\Delta_l^{(f)}$ Bus Angle at the From-End of Line $l$

$z_b^{(i)}$ Binary Decision Variable for Splitting Bus $i$, 1 Indicates the Bus Is Connected

$zb_l^{(f)}$ Binary Decision Variable for Splitting the From-End of Line $l$, 1 Indicates the Bus Is Connected

$z_{DG}^{(i)}$ Binary Decision Variable for the MCP Connection to Bus $i$, 1 Indicates It Is Connected to Side 1

$P_G^{(s)}$ Slack Bus Active Power Generation (MW)

$P_l^{l'}$ Active Power Flow of Line $l$ for the Contingency Scenario of Line $l'$ Removal

$P_{l_1}^{(f,l')}$ Active Power Flow of Line $l$ Connected to Side 1 of the From-End of the Line (MW) for the Contingency Scenario of Line $l'$ Removal

# I. Introduction

The integration of an increasing number of renewable energy sources into contemporary power grids has made Congestion Management (CM) crucial for the secure operation of these systems. Importantly, CM can incur significant costs [1], [2]; therefore, various methodologies have been explored and investigated to effectively address CM challenges and mitigate congestion. Two main approaches to addressing congestion in transmission networks involve either focusing on the root causes or employing congestion alleviation techniques, as can be seen in Fig. 1. Regarding the former, recent work by the present author in [1] seeks to shed light on the primary causes of congestion in transmission grids. The findings suggest that imbalances within the reactance structure of the grids are a significant factor contributing to congestion. However, these results are preliminary, and further investigation is needed to fully explore the underlying fundamental causes of congestion. Moreover, addressing the potential root cause of congestion does not guarantee complete relief, given the complex nature of power grids. Therefore, relying on congestion alleviation techniques remains essential for securely managing grid operations.

Regarding congestion alleviation, two primary strategies are employed to facilitate congestion relief. The first strategy is Network Topology Optimization (NTO) [3], which Transmission System Operators (TSO) may implement to mitigate the congestion. The NTO strategy may involve bus switching, transmission line switching, or, depending on the circumstances, a combination of both actions. Notably, when line switching is considered exclusively, it may also be referred to as Optimal Transmission Switching (OTS) [4]. Several methodologies have been presented in the literature for the NTO. For example, [5] introduces a novel topology control algorithm designed to alleviate congestion and overloads in power networks, thereby enhancing security margins. The algorithm employs a heuristic based on the normalized squares of line currents as an indicator of transmission line loading. By assessing the sensitivity of this measure relative to breaker positions, the heuristic identifies optimal switching actions to reconfigure the network topology, reducing congestion and overloads while maintaining a largely fixed generation dispatch. [6] proposes a topology-aware data-driven Optimal Power Flow (OPF) framework capable of adapting to network topology changes without requiring extensive model retraining. As one of the main constraints of the proposed optimization framework, the transmission line loading must satisfy the corresponding maximum loading capacity limits. The proposed methodology consists of two main stages. In the first stage, a topology transfer framework is developed to generate sample-efficient data-driven OPF models under new network topologies. Gaussian Process Regression (GPR) is utilized to capture the nonlinear relationship between the actual OPF solutions obtained from conventional optimization solvers and the predictions produced by previously trained data-driven OPF models.

In the second stage, a dynamic ensemble learning strategy is introduced, in which both the ensemble weights and the selected topology transfer models are adaptively determined in order to improve prediction accuracy and avoid biased OPF solutions. Another work is [7] where the challenges faced by TSO in ensuring secure and cost-effective power system operations, particularly in the presence of transmission congestion, which can hinder the dispatch of more economical generating units are investigated. The study introduces a decentralized transmission switching methodology for congestion management in interconnected multi-area power systems. The proposed framework coordinates switching actions across both regional networks and tie lines while accounting for credible contingencies, enabling each regional control center to independently manage congestion within its operational area. The performance of the method is demonstrated using a range of benchmark systems..

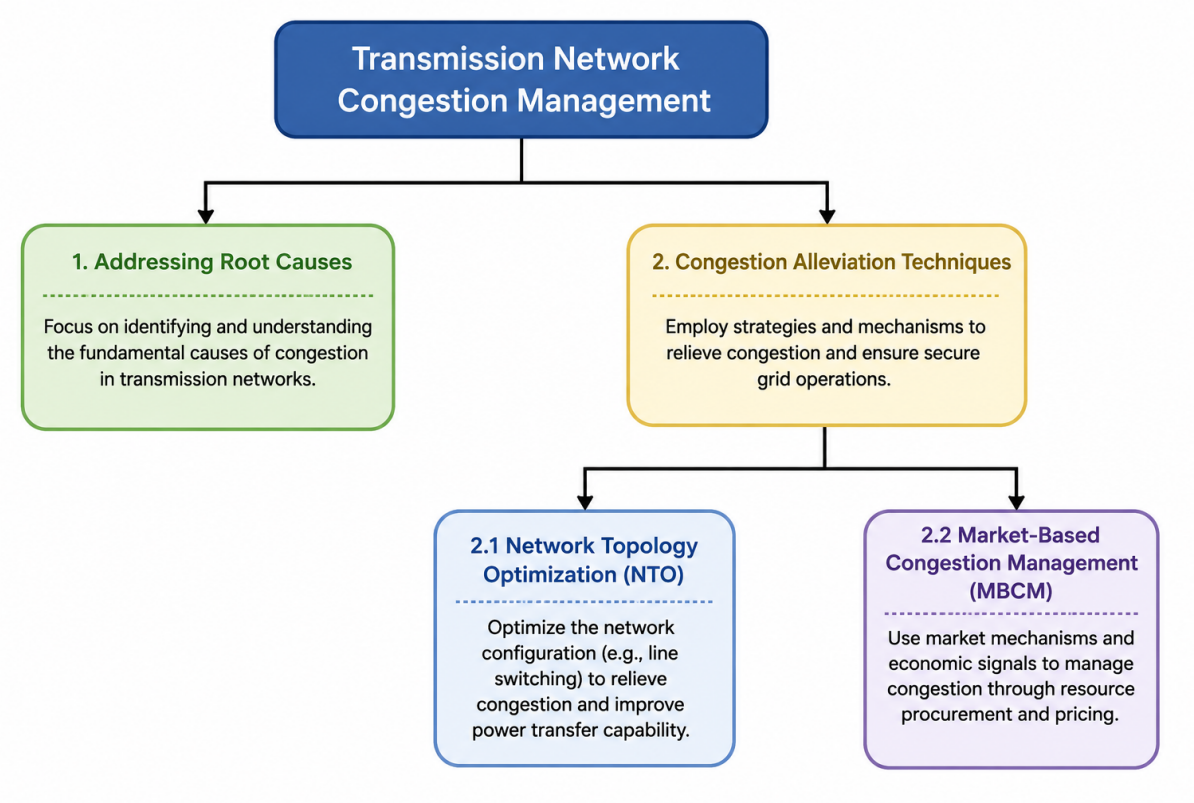

Fig. 1. Categorization of Transmission Network Congestion Management Approaches

The work in [3] proposes an advanced NTO model, based on substation and node-breaker representations, to effectively mitigate congestion and reduce operational costs in power systems. Unlike the traditional bus-branch model, which lacks the ability to fully capture detailed substation switching actions (such as bus splitting), the proposed NTO model offers a more practical and realistic depiction of these operations, providing greater flexibility in switching actions. [8] introduces a heuristic approach combining optimal transmission line switching and bus splitting to alleviate congestion and reduce operational costs by rerouting power flows within the network. The work presented in [8] develops a topology optimization framework that combines bus splitting and transmission line switching to alleviate network congestion while maintaining system security. The proposed methodology formulates the problem as a two-level optimization model, where the upper level determines the optimal switching actions using a mixed-integer second-order cone programming formulation, and the lower level evaluates the feasibility of the resulting network topology under $N-1$ contingency conditions. The study further highlights that, despite the effectiveness of bus splitting and its relatively limited operational impact compared with line switching, its practical application has been constrained by the complexity associated with detailed node-breaker modeling.

The second congestion alleviation strategy involves employing a congestion management market mechanism, also referred to as Market Based Congestion management (MBCM). In this approach, power market participants modify their operational settings in accordance with a predetermined market framework. This adjustment facilitates operators in effectively managing and alleviating congestion within the power system. The congestion management market mechanism is generally designed to support an unconstrained power market framework, where the initial market clearing process is performed without explicitly incorporating transmission network operational constraints, similar to the practice adopted in many European electricity markets. This approach addresses the potential congestion that may arise from the participants' initial accepted operational setpoints. Various congestion management market mechanisms have been developed in the literature. For example, [9] presents a novel platform that integrates P2P trading with secure, automated low-voltage grid control to manage congestion effectively. The proposed dual local energy market includes a P2P energy market and a flexibility market, the latter of which proactively addresses congestion. Additionally, the platform offers real-time, curative congestion management and enables prosumers, users who both produce and consume energy, to optimize their financial returns by engaging in both markets. Another example is [10] in which the authors examine market-based CM in power systems, focusing on enabling flexibility from both small and large-scale resources through a platform that integrates this flexibility into markets for use by system operators (Distribution System Operators (DSO) and TSO). The study identifies two key challenges: low market liquidity and the adverse effects of flexibility activation. To address these issues, the paper suggests enhancing market integration and DSO-TSO coordination within the platform. By incorporating bids from the Nord Pool intraday and balancing markets, the study demonstrates that established markets can help alleviate liquidity constraints. [11] introduces a stochastic model utilizing chance constraints for managing grid congestion in the day-ahead power market. By accounting for the uncertainties in wind generation and demand-side response, the model optimizes daily generator dispatch and demand-responsive loads to minimize the risk of grid congestion. The stochastic model is converted into a deterministic equivalent by leveraging the probability distribution of random variables, while also considering transmission system losses. The probabilistic power flow is used to assess the satisfaction of chance constraints and quantify congestion risk via the congestion risk value. Simulations of the test systems validate the model's effectiveness for congestion management in the day-ahead market.

It is evident that the foundation of all CM market mechanisms is the flexibility that participants can offer and trade within the market, enabling the TSO to effectively manage congestion. A substantial source of potential flexibility in modern power grids resides in the capabilities offered by aggregators [2], [12]. Through the operation and management of a diverse array of loads, generators, and Energy Storage Systems (ESS), aggregators demonstrate significant flexibility, positioning them as crucial participants in effective congestion management market mechanisms. Consequently, any MBCM methodology seems to need to incorporate the potential flexibility provided by aggregators. Hence, they have garnered considerable attention within contemporary CM strategies for today's power grids. Notably, aggregators can be categorized based on the number of their access points to transmission lines, distinguishing between those with a Single Connection Point (SCP) and those with Multiple Connection Points (MCP) [2]. In this regard, a significant pioneering work investigating the flexibility services provided by aggregators for CM is presented in [2], where the authors discuss both SCP and MBCM in electric power grids using aggregators of dispersed small-scale consumers, producers, and prosumers within the distribution network. The proposed method utilizes a market-driven framework where the grid operator conducts minimum-cost re-dispatching of generators, loads, and aggregators to alleviate congested lines. The power network is modeled using a linearized DC approach, applying Power Transfer Distribution Factors (PTDFS) and Topological Generation and Load Distribution Factors (TGDFS and TLDFS).

From the aforementioned discussion, two primary sources of potential flexibility to address the congestion management problem emerge: the flexibility offered by the NTO and that provided by power market participants, including aggregators.

The main contribution of this paper is the development of an innovative congestion management framework that simultaneously exploits the operational flexibility provided by NTO and the market-based flexibility services offered by various participants, particularly aggregators. To achieve this, a novel Integrated Congestion Management Power Market (ICMPM) is proposed, in which market participants are allowed to submit flexibility offers while the TSO can utilize NTO actions to minimize the overall CM cost. The proposed framework is formulated using MILP and explicitly incorporates the $N-1$ security constraints, fundamental constraints to ensure secure grid operation under contingency conditions. Furthermore, aggregators are modeled according to their access points to the transmission network, enabling the coordinated integration of different categories of flexibility services within the proposed market framework. In this manner, the proposed ICMPM effectively bridges the operational flexibility of NTO with the substantial flexibility resources available in the electricity market. Notably, many existing MBCM methodologies, such as the approach presented in [2], rely on topological sensitivity matrices including PTDF, LODF, GGDF, and GLSK. Such methods face significant challenges in incorporating NTO into the market framework, since the admittance structure of the grid becomes variable due to topology changes, thereby rendering matrix-based methods unsuitable for such a dynamic environment. In contrast, the formulation proposed in this paper effectively overcomes this limitation and enables the successful integration of NTO within the congestion management framework.

## II. Methodology

### A. The Integrated Congestion Management Power Market

In this section, the novel Integrated Congestion Management Power Market (ICMPM) is introduced. As previously discussed, this market is underpinned by two primary forms of flexibility: flexibility provided by power market participants, particularly aggregators, and that offered through NTO. In the following, the two fundamental flexibilities inherent in the proposed ICMPM are first examined. Subsequently, the formulation and structural configuration of the market are articulated and discussed.

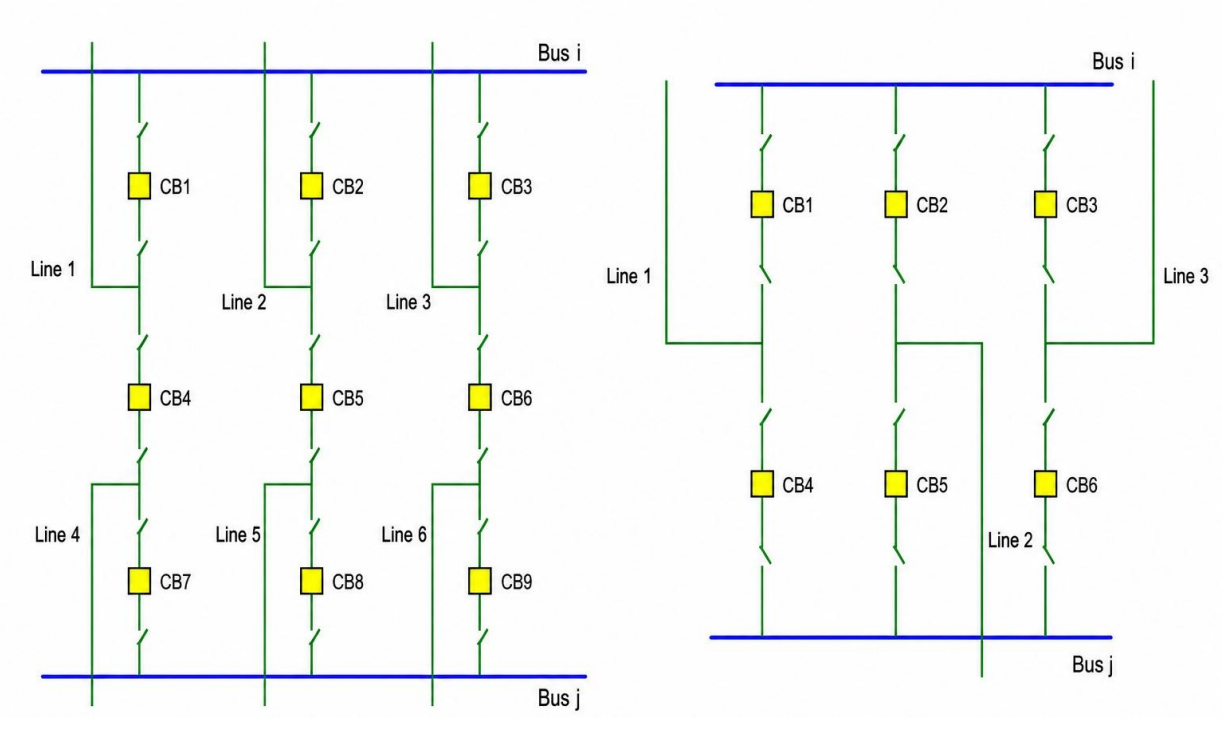

Fig. 2. Substation bus configurations; (a) Breaker and Half Configuration, (b) Double Bus Double Breaker Configuraution

*1) The flexibility conferred by the aggregators:* While the flexibility offered by conventional market participants, such as traditional generators, has been extensively examined in existing literature, the emerging and significant flexibility potential of aggregators warrants comprehensive modeling. Detailed modeling of this flexibility is essential for the development of new congestion management markets.

Aggregators typically manage and operate a diverse range of consumers, prosumers, and small generating units within distribution networks. By aggregating the capacity of these resources, they can offer flexibility to the TSO at their respective connection points to the transmission networks. The work presented in [2] offers a comprehensive categorization of aggregators based on their specific features, a framework that will also be adopted and considered in the present study. In this sense, the aggregators are categorized into two main types: the aggregators with Single Connection Point (SCP) to the transmission grids and the aggregators with Multiple Connection Points (MCP) to the transmission grids. Aggregators submit their flexibility to the TSO at specific prices, as their offers in the congestion management power market, determined through market mechanisms at distribution grids such as the one presented in [2] when they aggregate their resources. The internal mechanism for determining the offers of these aggregators falls outside the scope of this paper; it is assumed that the aggregators have already prepared their offers.

Moreover, depending on the power market regulations, participants in the congestion market may submit different types of offers, also known as "block orders". Typically, four main conventional types facilitated in power markets are: "all-or-nothing" block orders, linked block orders, curtailable block orders, and profile block orders [2]. Notably, different power markets may facilitate various types of block orders, but two that are widely used across multiple power markets are all-or-nothing and curtailable block orders, which are therefore considered in this paper. For all-or-nothing blocks, the entire block order must either be fully accepted or rejected by the market operator. In contrast, for curtailable blocks, the market operator can curtail and partially accept the block orders within the limitations specified by the participants.

Regarding the block orders for aggregators, as discussed in [2], a significant challenge arises for MCP aggregators due to their requirement to either fully accept or reject a specified amount of power or load change at various locations. Consequently, [2] proposes the use of "location-related block orders" for MCP aggregators, which incorporate an "all-or-nothing" feature, thereby ensuring that the flexibility offered is specific to the corresponding locations. The present paper similarly examines location-oriented block orders for MCP aggregators. Notably, SCP aggregators can be treated similarly to conventional generating units, allowing them to submit either all-or-nothing or curtailable offers.

*2) The flexibility conferred by NTO:* The NTO may involve both bus switching and line switching. Typically, NTO is

conducted at substations, where appropriate switching at busbars enables the desired configuration. As shown in Fig. 2, the two common arrangements for busbars at a substation are the double-bus double-breaker and the breaker-and-a-half configurations [3], [8]. Notably, while the double-bus double-breaker configuration offers greater flexibility, it is more expensive than the breaker-and-a-half arrangement due to the increased number of breakers required. As a result, the breaker-and-a-half configuration is more commonly used in power grids.

*3) The structure and formulation of the proposed* ICMPM*:* The principal challenge lies in integrating the NTO capability with MCP and SCP aggregators, as well as other conventional market participants, to effectively develop an integrated congestion management power market. The conventional approach for facilitating the congestion market, even the ones with aggregators like the work in [2], involves determining the load or generation changes for each bus and then analyzing the effects of these changes on grid lines using matrices like PTDF, LODF, GGDF, GLSK, etc to derive an optimal solution. However, when incorporating the NTO into the market, the admittance structure of the grid becomes variable due to changes in topology, rendering the matrix-based methods unsuitable for this dynamic context.

Conversely, although the NTO problem has already been formulated utilizing MILP in [3], this conventional NTO model fails to incorporate the flexibility offered by aggregators as a significant source of capacity within modern power grids. This limitation hampers the model's ability to adequately address congestion management in power markets. Furthermore, it neglects to integrate essential operational security constraints, notably the $N-1$ security constraints, which is typically mandated by the TSO to ensure system robustness and reliability under contingencies. Building upon the MILP-based NTO modeling approach presented in [3], the methodology developed in this paper extends the formulation to effectively address the aforementioned requirements within the proposed congestion management framework.

Regarding the power market design for integrating the proposed ICMPM into existing electricity market frameworks, one may refer to the structural model illustrated in Fig. 3. This model is based on the generalized and simplified unconstrained power market architecture commonly implemented in European countries [13]. While the proposed structure aligns with the core characteristics of most European power markets and those adhering to similar principles, variations in implementation may occur from country to country. Consequently, adaptations of the proposed market can be tailored by market participants according to specific national requirements and operational conditions. For instance, in power markets such as the Irish market's I-SEM, where the current 'balancing market' possesses sufficient capacity to effectively address congestion events, the proposed ICMPM can be incorporated within the existing balancing market infrastructure (like the one in Fig. 3), thereby

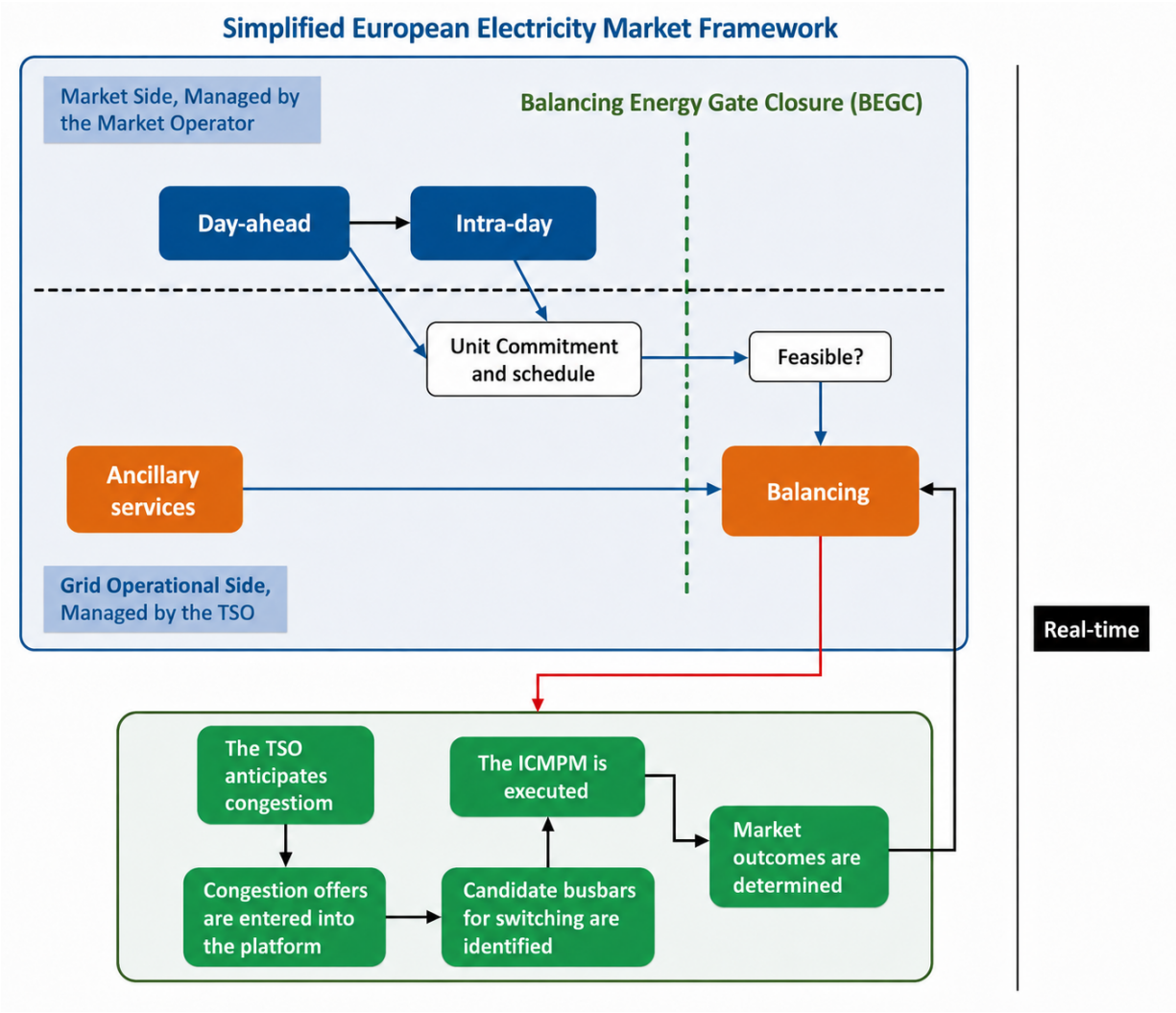

Fig. 3. Integration of the proposed ICMPM scheme into the European electricity market setup

enhancing operational efficiency without necessitating extensive structural overhauls.

Regarding the temporal feasibility of the proposed ICMPM, two principal aspects of the methodology should be considered. First, with respect to the flexibility services provided by aggregators, it has been demonstrated in [2] that a time interval of approximately one hour prior to real-time operation, corresponding to the typical BEGC timeframe adopted in markets such as I-SEM, is sufficiently practical for aggregators to submit their offers, for the market-clearing process to be executed, and subsequently for the corresponding BM dispatch notifications to be communicated. Second, concerning the implementation timeframe of topological remedial actions, NTO, it is shown in [14] that such actions can be executed within sufficiently short timescales. Depending on the operational characteristics and control infrastructure of the transmission network, the activation time may range from several minutes up to approximately one hour. Therefore, the one-hour BEGC timeframe adopted in several European electricity markets, such as the I-SEM, is well aligned with the operational requirements of the proposed framework. Consequently, the integrated mechanism of the proposed ICMPM satisfies the temporal and operational requirements of European electricity markets, as well as other market structures operating under comparable scheduling and balancing timeframes.

To enable NTO capability, the configuration presented in Fig.4 for each candidate bus is considered. This adopted bus configuration facilitates the connection of each power market participant, including aggregators, to either side 1 or side 2 of the respective buses. This arrangement ensures same connectivity options for all regular loads and allows for conventional line-switching and bus-splitting operations. The market players considered for the proposed integrated congestion management market include MCP and SCP aggregators, as well as conventional generating units and loads.

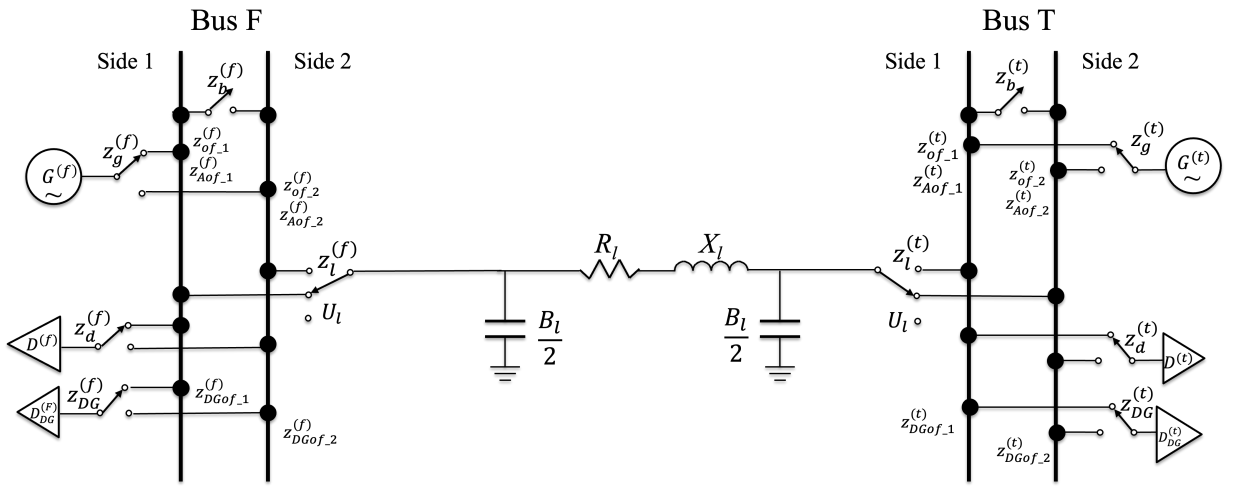

Fig. 4. Schematic Representation of the Adopted Bus Switching Configuration With Aggregator Flexibility Services Enabled

The formulation of the proposed novel integrated congestion management power market is developed using MILP programming and is presented as follows:

$$
\begin{aligned}
&Min\Big[\sum_{m}((P^{(m)}_{Gof_1}+P^{(m)}_{Gof_2})\times C^{(m)}_{Gof})+\\
&\sum_{k}((z^{(k)}_{Aof_1}+z^{(k)}_{Aof_2})\times C^{(k)}_{Aof})+\\
&\sum_{w}((z^{(w)}_{DGof_1}+z^{(w)}_{DGof_2})\times C^{(w)}_{DGof})+M\times P^{(s)}_{G}\Big]\\
&\quad, \forall m\in G^{of}\ ,\ \forall k\in G^{Aof}\ ,\ \forall w\in D^{DG}
\end{aligned}
\tag{1}
$$

$$
\begin{aligned}
&\sum\big[P^{(i)}_{l_1}-(z^{(i)}_{g}\times G^{(i)})-P^{(i)}_{Gof_1}\\
&-(z^{(i)}_{Aof_1}\times P^{(i)}_{Aof})+(z^{(i)}_{d}\times D^{(i)})+\\
&(z^{(i)}_{DGof_1}\times P^{(i)}_{DGof})\big]=0\,,\\
&\quad, \forall i\in B\ ,\ \forall l\in L
\end{aligned}
\tag{2}
$$

$$
-z^{(i)}_{of_1}\times G^{(i)}_{max}\le P^{(i)}_{Gof_1}\le z^{(i)}_{of_1}\times G^{(i)}_{max}
\tag{3}
$$

$$
\begin{aligned}
&\sum\big[P^{(i)}_{l_2}-((1-z^{(i)}_{g})\times G^{(i)})-P^{(i)}_{Gof_2}\\
&-(z^{(i)}_{Aof_2}\times P^{(i)}_{Aof})+((1-z^{(i)}_{d})\times D^{(i)})+\\
&(z^{(i)}_{DGof_2}\times P^{(i)}_{DGof})\big]=0\,,\\
&\quad, \forall i\in B\ ,\ \forall l\in L
\end{aligned}
\tag{4}
$$

$$
-z^{(i)}_{of_2}\times G^{(i)}_{max}\le P^{(i)}_{Gof_2}\le z^{(i)}_{of_2}\times G^{(i)}_{max}
\tag{5}
$$

$$
\begin{aligned}
&-z^{(f)}_{l}\times P^{max}_{l}\ \le P^{(f)}_{l_1}\le z^{(f)}_{l}\times P^{max}_{l}\\
&-(1-z^{(f)}_{l})\times P^{max}_{l}\ \le P^{(f)}_{l_2}\le(1-z^{(f)}_{l})\times P^{max}_{l}\\
&, \forall l\in L
\end{aligned}
\tag{6}
$$

$$
\begin{aligned}
&-z^{(t)}_{l}\times P^{max}_{l}\ \le P^{(t)}_{l_1}\le z^{(t)}_{l}\times P^{max}_{l}\\
&-(1-z^{(t)}_{l})\times P^{max}_{l}\ \le P^{(t)}_{l_2}\le(1-z^{(t)}_{l})\times P^{max}_{l}\\
&, \forall l\in L
\end{aligned}
\tag{7}
$$

$$
\begin{aligned}
&P_{l}=P^{(f)}_{l_1}+P^{(f)}_{l_2}\\
&P_{l}=P^{(t)}_{l_1}+P^{(t)}_{l_2}\\
&P_{l}\times X_{l}=\Delta^{(f)}_{l}-\Delta^{(t)}_{l}\\
&, \forall l\in L
\end{aligned}
\tag{8}
$$

$$
\begin{aligned}
&-(1-U_{l})\times P^{max}_{l}\le P_{l}\le(1-U_{l})\times P^{max}_{l}\\
&1-U_{l}-z^{(f)}_{l}\ge 0\\
&1-U_{l}-z^{(t)}_{l}\ge 0\\
&, \forall l\in L
\end{aligned}
\tag{9}
$$

$$
\begin{aligned}
&-(1-z^{(f)}_{l})\times\delta_{max}\le\Delta^{(f)}_{l}-\delta^{(f)}_{1}\le(1-z^{(f)}_{l})\times\delta_{max}\\
&-z^{(f)}_{l}\times\delta_{max}\le\Delta^{(f)}_{l}-\delta^{(f)}_{2}\le z^{(f)}_{l}\times\delta_{max}\\
&-(1-z^{(t)}_{l})\times\delta_{max}\le\Delta^{(t)}_{l}-\delta^{(t)}_{1}\le(1-z^{(t)}_{l})\times\delta_{max}\\
&-z^{(t)}_{l}\times\delta_{max}\le\Delta^{(t)}_{l}-\delta^{(t)}_{2}\le z^{(t)}_{l}\times\delta_{max}\\
&, \forall l\in L
\end{aligned}
\tag{10}
$$

$$
\begin{aligned}
&-\delta_{max}\times(1-z^{(i)}_{b})\le\delta^{(i)}_{2}-\delta^{(i)}_{1}\le\delta_{max}\times(1-z^{(i)}_{b})\\
&, \forall i\in B
\end{aligned}
\tag{11}
$$

$$
\begin{aligned}
&zb^{(f)}_{l}\le z^{(f)}_{l}\\
&zb^{(t)}_{l}\le z^{(t)}_{l}\\
&, \forall l\in L
\end{aligned}
\tag{12}
$$

$$z_b^{(i)} \leq z_g^{(i)}$$
$$z_b^{(i)} \leq z_d^{(i)}$$
$$z_b^{(i)} \leq z_{DG}^{(i)}$$
$$z_g^{(i)} \geq z_{of_1}^{(i)}$$
$$(1 - z_g^{(i)}) \geq z_{of_2}^{(i)}$$
$$z_{DG}^{(i)} \geq z_{DGof_1}^{(i)} \quad (13)$$
$$(1 - z_{DG}^{(i)}) \geq z_{DGof_2}^{(i)}$$
$$z_g^{(i)} \geq z_{Aof_1}^{(i)}$$
$$(1 - z_g^{(i)}) \geq z_{Aof_2}^{(i)}$$
$$, \forall i \in B, \ i \neq s$$

$$\sum (P_{l_1}^{(s)}) - P_G^{(s)} + D^{(s)} = 0 \,, \forall l \in L$$
$$zb^{(s)} = 1 \quad (14)$$
$$\delta_1^{(s)} = 0$$

$$\forall l' \in L' :$$
$$P_{l_n}^{(i)} \to P_{l_n}^{(i,l')}$$
$$P_{l_n}^{(t)} \to P_{l_n}^{(t,l')}$$
$$P_{l_n}^{(f)} \to P_{l_n}^{(f,l')}$$
$$P_l \to P_l^{(l')}$$
$$\delta_n^{(i)} \to \delta_n^{(i,l')} \quad (15)$$
$$\Delta_l^{(f)} \to \Delta_l^{(f,l')}$$
$$\Delta_l^{(t)} \to \Delta_l^{(t,l')}$$
$$P_G^{(s)} \to P_G^{(s,l')}$$

*Regenerate constraints EQ2, EQ4, EQ6 – EQ12, EQ14*

$$P_l^{(l')} = P_{l_n}^{(t,l')} = P_{l_n}^{(f,l')} = 0 \ \forall l' \in L', n \in \{1, 2\} \quad (16)$$

It is assumed that, in the formulation, $P_{Gof_1}^{(i)}$, $P_{Gof_2}^{(i)}$, $P_{DGof}^{(i)}$, $P_{Aof}^{(i)}$, $G^{(i)}$, and $D^{(i)}$ are bus-indexed vectors whose dimensions are equal to the total number of buses. For buses at which the corresponding entity is not present, the associated entry is set to zero; otherwise, the entry takes the corresponding value associated with that bus. Equation (1) represents the objective function of the proposed model, which aims to minimize the total congestion market costs. The constraints defined in equations (2) to (5) correspond to the nodal injection constraints for both sides (side 1 and side 2) of all buses. Equations (6) and (7) establish the criteria by which each line is connected to either side 1 or side 2 at the respective *From-end* and *To-end*. The constraints outlined in equation (8) ensure that the power magnitude at the *From-end* and *To-end* corresponds to the line power flow, which is also computed based on the DC power flow assumptions. Equation (9) introduces the line switching constraints. The constraints in equation (10) determine the bus angle at both the *From-end* and *To-end*, by identifying whether each corresponding bus is associated with side 1 or side 2. Equations (11) and (12) enable the connection or disconnection of each bus. The constraints outlined in Equation (13) facilitate the allocation of the flexibility services provided by aggregators, namely curtailable generation for SPC aggregators (or conventional generators) and all-or-nothing generation adjustment for both SCP and MCP aggregators, to the corresponding buses, while simultaneously enabling bus-switching operations. Equation (14) pertains to the slack bus, and finally, equations (15) and (16) impose the $N-1$ contingency constraints on the formulation. Notably, the implementation of the $N-1$ constraints necessitates the modification of certain variable notations to facilitate iteration across all considered $N-1$ contingency scenarios, as outlined in the equation (15).

Moreover, across all $N-1$ contingency scenarios, the decision variables associated with both the generation flexibilities provided by aggregators and the bus-switching actions remain identical. This enables the determination of a unified set of decision variables for congestion relief while ensuring that the $N-1$ security constraints are satisfied. Regarding the $N-1$ contingency analysis, it is generally neither necessary nor computationally efficient to consider all possible contingency scenarios, since the incorporation of additional $N-1$ security constraints typically leads to more conservative operational decisions and consequently higher congestion management costs [1]. In practical power system operation, TSO generally performs contingency selection procedures in advance, based on operational practices, reliability criteria, and grid-code requirements, in order to identify the set of contingencies against which the system operation must remain secure. Consequently, different subsets of $N-1$ contingency scenarios may be adopted depending on the operational requirements and characteristics of the considered transmission network.

In this paper, all non-radial transmission lines are considered as potential $N-1$ contingency scenarios. Radial lines are excluded from the contingency set since their outage does not result in power-flow redistribution capable of causing overloads elsewhere in the network and, therefore, they do not significantly influence congestion management requirements [1]. To identify the radial lines, the methodology based on the *radiality index*, as presented in [1], [15], is employed. This methodology utilizes the PTDF matrix, $\Gamma$, as follows:

$$\begin{aligned} &\Omega(\ell\,,\,:) = \Gamma(:\,\,,\,j) \\ &\text{S.t.:} \\ &\quad n_r = i, \quad P_i = 1, \quad P_j = -1 \end{aligned} \tag{17}$$

Where $\Omega$ is the *Line Mesh Matrix* [15]. Note that in the $\Omega$ matrix, for all $m, n \in B$ and $m \neq n$, if $\Omega(n,\ m) \leq \Omega(m,\ n)$, then $\Omega(n,\ m)$ is replaced with $\Omega(m,\ n)$. The diagonal components in $\Omega$ represent the *radiality indices* of the corresponding lines, where a diagonal component equal to one indicates that the corresponding line is radial.

$$\begin{aligned} &\mathbb{R}_l = \Omega(\ell\,,\,\ell) \\ &0 \leqslant \mathbb{R}_l \leqslant 1 \end{aligned} \tag{18}$$

Finally, it should be noted that, depending on practical and operational requirements, additional constraints may be incorporated into the proposed model, and supplementary terms may be included in the objective function. For instance, if an operational cost is associated with switching actions, the following term may be added to the objective function presented in Equation (1):

$$\begin{aligned} &\sum_h \left(C_{sw}^{(h)} \times z_b^{(h)}\right) \\ &, \forall h \in B^{sw} \end{aligned} \tag{19}$$

Furthermore, constraints presented in equations (20) and (21) respectively limit the number of line-switching and bus-switching operations to predefined allowable values:

$$\begin{aligned} &\sum_l \left(U_l\right) \leq F_l \\ &, \forall l \in L \end{aligned} \tag{20}$$

$$\begin{aligned} &\sum_i \left(1 - Z_b^{(i)}\right) \leq F_b \\ &, \forall i \in B \end{aligned} \tag{21}$$

Another potential modification arises when the TSO chooses not to permit switching actions for a specific set of buses, or when bus switching is not feasible due to operational or technical limitations. In such cases, the following constraints can be incorporated into the proposed model:

$$Z_b^{(j)} = 1\,,\ \forall j \in B_{sw}^{'} \tag{22}$$

## III. Results

To implement the methodology presented in this paper, two toolboxes within the MATLAB environment [16] are utilized: 1) The handling and augmentation of the test system data are performed using the MATPOWER toolbox [17]; 2) Furthermore, the novel ICMPM presented in SectionII-A3 is modeled using YALMIP [18]. The underlying script developed for implementing the proposed ICMPM framework is publicly available through the online repository referenced in [19].

### A. Test platform

From the repository presented in [20], the `pglib_opf_case39_epri` test system is selected as the case study for the assessment of the proposed framework. The original loading profile of the grid prior to the implementation of any congestion management action can be observed in Fig. 5. As can be seen from this figure, the grid experiences congestion on several transmission lines. Also, the flexibility services provided by aggregators through block orders are summarized in Table. I. It should be noted that the offer values considered in this study are hypothetical and are introduced solely to facilitate the evaluation and validation of the proposed ICMPM framework. Therefore, these values do not represent actual market offers observed in real-world electricity markets. For practical implementations, the proposed framework can readily incorporate real market bids and offers obtained from operational market environments.

Regarding the flexibility services provided by aggregators in Table I, the MCP offer effectively represents a load reduction action, as commonly encountered in demand-side management strategies, which may equivalently be interpreted as a reduction in electricity generation. As indicated in the table, it is assumed that this block order follows an all-or-nothing structure, meaning that, if accepted, a 200 MW load reduction is simultaneously activated at each of the specified buses; otherwise, no corrective action is undertaken. In addition, a curtailable SCP aggregators (or equivalently a conventional dispatchable generators) are represented by Agr_3 to Agr_10, with a maximum generation adjustments capability of ±70 MW. This implies that each generator or aggregator belonging to this offer category is connected to its corresponding bus and is capable of performing generation adjustments within its specified operating limits.

Regarding the operational flexibility associated with NTO, it is assumed that the maximum number of allowable switching actions is limited to four. No restrictions are imposed on the selection of candidate buses. However, as discussed in Section II-A3, using Equation 22 one may restrict the set of candidate buses according to their specific operational requirements and preferences of the TSO.

TABLE I
FLEXIBILITY OFFERS PROVIDED BY AGGREGATORS

| Offer ID | Connected Bus | Flexibility (MW) | Price (€/MWh) | Type |
|---|---|---|---|---|
| Agr_1 | [16,26,39] | -200 | 8 | MCP all-or-nothing |
| Agr_2 | 30 | +30 | 15 | SCP all-or-nothing |
| Agr_3 - Agr_10 | [32 - 39] | ± 70 | [20, 23, 24, 27, 31, 35, 21, 44] | SCP curtailable |

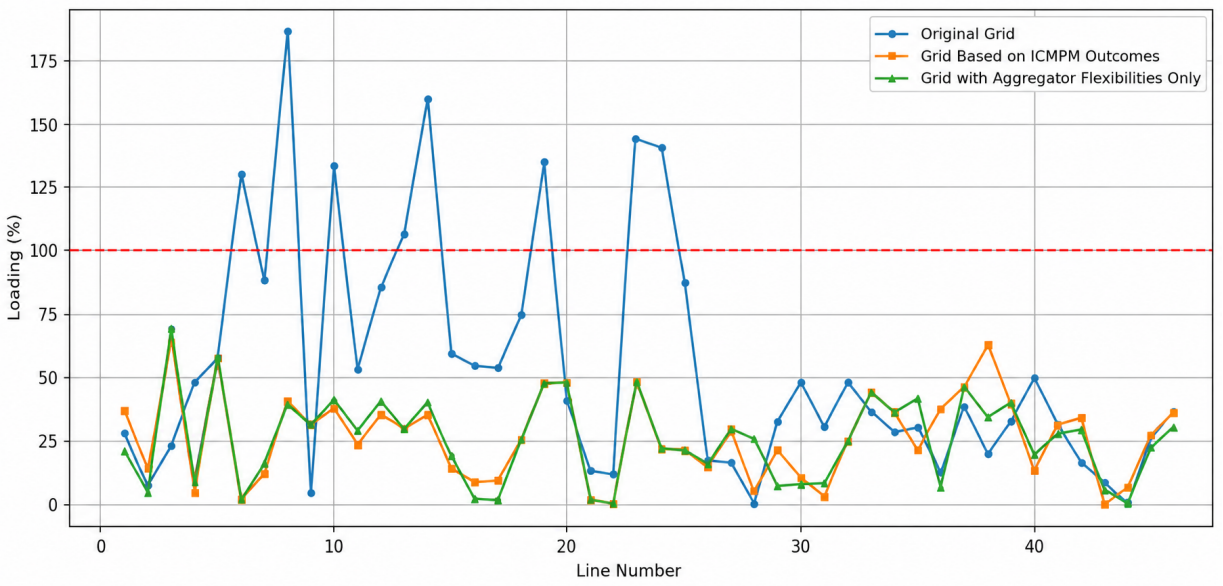

Fig. 5. Comparison of Line Loading Profiles for the `pglib_opf_case39_epri` Test System Under Different Congestion Management Scenarios

TABLE II
ACCEPTED OFFERS PROVIDED BY AGGREGATORS

| Offer ID | Accepted Bus Number | Accepted Flexibility (MW) | Bus Side (if applicable) |
|---|---|---|---|
| Agr_1 | [16, 26, 39] | -200 | [1, N/A, N/A] |
| Agr_2 | Not Accepted | N/A | N/A |
| Agr_7, Agr_9, Agr_10 | [36, 38, 39] | [60, -2.83, 70] | [N/A, N/A, N/A] |

The outcomes of the implementation of the proposed ICMPM are illustrated in Fig.6, in which the buses denoted in red correspond to the buses selected for bus-switching operations. Note that Bus31, highlighted in green, represents the slack bus of the system.

The details of the accepted offers obtained through the proposed ICMPM are presented in Table II. It should be noted that, in this table, if an offer associated with a particular bus is accepted while the corresponding entry in the 'Bus Side' column is indicated as 'N/A', this implies that although the aggregator's offer has been accepted, no bus-switching operation has been performed at that bus. Consequently, the bus-side designation is not applicable in such cases.

Notably, as an example, Fig. 7 illustrates the configuration of Bus 16 based on the outcomes of the proposed ICMPM.

From the market outcomes presented in Table. II, it can be observed that the MCP offer Agr_1 is accepted, whereas the all-or-nothing offer Agr_2 is not accepted. In addition, three curtailable SPC offers associated with Agr_7, Agr_9 and Agr_10 are also accepted. Notably, as can be observed in Fig. 5, the proposed ICMPM successfully relieves the congestion throughout the network.

To investigate the effects of integrating NTO within the proposed congestion management market, the proposed ICMPM is executed again, this time by enforcing the constraint $\sum z_b^{(i)} = N^B$, where $N^B$ is the number of buses, thereby prohibiting any switching actions. Table. III presents the resulting market outcomes based on the same set of received offers provided in

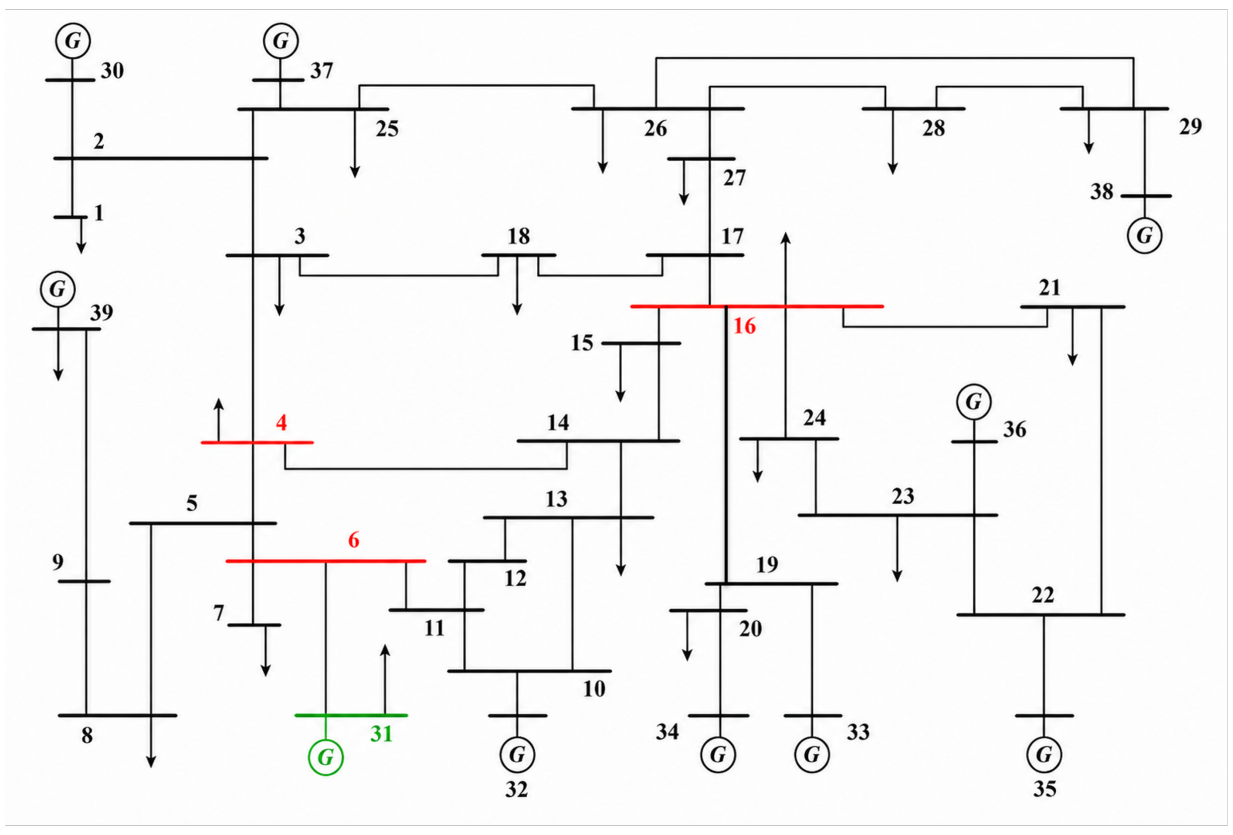

Fig. 6. `pglib_opf_case39_epri` test system with red buses selected for bus switching and the green bus representing the slack bus.

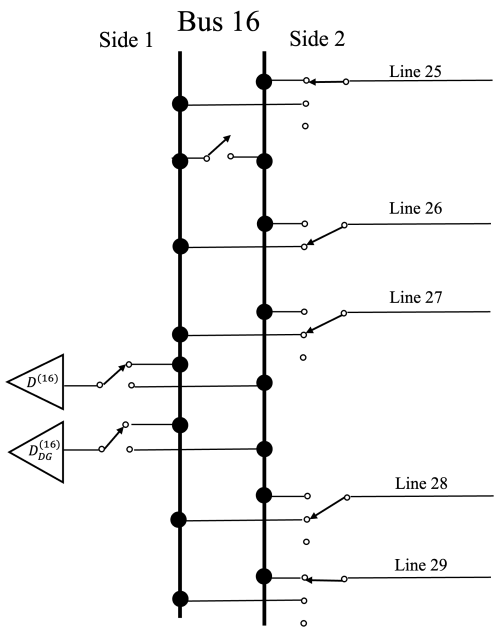

Fig. 7. Configuration of Bus 16 Based on the Outcomes of the Proposed ICMPM

TABLE III
ACCEPTED AGGREGATOR OFFERS WITHOUT ALLOWING SWITCHING OPERATIONS

| Offer ID | Accepted Bus Number | Accepted Flexibility (MW) |
|---|---|---|
| Agr_1 | [16, 26, 39] | -200 |
| Agr_2 | Not Accepted | N/A |
| Agr_3 - Agr_10 | [32 - 39] | [70, 70, 70, 70, 70, -35.18, -70, 70] |

Table. I. As can be observed in Fig. 5, the flexibility services provided by the aggregators successfully relieve the congestion within the grid.

Significantly, the total flexibility cost in the congestion management market when allowing NTO is reported to be 11,605€, whereas the corresponding cost without allowing NTO is reported to be 14,560€, representing approximately a 20.3% reduction in the overall congestion management cost. This clearly demonstrates the benefits of integrating the operational flexibility provided by NTO with the flexibility services offered by aggregators.

## IV. Conclusions

This paper presented a novel Integrated Congestion Management Power Market (ICMPM) framework that facilitates

the integration of flexibility services provided by aggregators, which constitute a significant characteristic of modern power systems, while simultaneously utilizing the operational flexibilities available to the TSO, namely NTO, under critical $N-1$ security constraints. The proposed ICMPM framework was formulated using MILP with the objective of minimizing the total flexibility procurement cost, while also providing the capability to incorporate additional constraints in order to satisfy specific operational requirements of a given power network. Furthermore, two categories of aggregator offer connections, namely Multiple Connection Point and Single Connection Point offers, together with two categories of order blocks, including all-or-nothing and curtailable block orders, were considered in the proposed framework. The proposed congestion management framework was evaluated using the `pglib_opf_case39_epri` test system as the case study, demonstrating the effectiveness of the proposed ICMPM in alleviating congestion scenarios. In addition, further analysis was conducted to investigate the impact of integrating NTO within the flexibility services provided by aggregators. The obtained results demonstrated that the incorporation of NTO can reduce the total flexibility procurement cost in the congestion management market by approximately 20%, thereby highlighting the economic and operational benefits of combining network operational flexibility with market-based flexibility services.

## References

[1] A. Beiranvand and P. Cuffe, "A congestion alleviation technique exploiting structural insights on the interaction of line loading limits, reactance and line outage security constraints," *IEEE Transactions on Power Systems*, vol. 38, no. 4, pp. 3719–3732, 2023.

[2] a. . T. MiloInternational Journal of Panto, keywords = Aggregator, Linear programming, Market-based congestion management, Optimization, Prosumer, "Market-based congestion management in electric power systems with exploitation of aggregators," *International Journal of Electrical Power & Energy Systems*, vol. 121, p. 106101, 2020.

[3] M. Heidarifar and H. Ghasemi, "A network topology optimization model based on substation and node-breaker modeling," *IEEE Transactions on Power Systems*, vol. 31, no. 1, pp. 247–255, 2016.

[4] M. Numan, M. F. Abbas, M. Yousif, S. S. M. Ghoneim, A. Mohammad, and A. Noorwali, "The role of optimal transmission switching in enhancing grid flexibility: A review," *IEEE Access*, vol. 11, pp. 32 437–32 463, 2023.

[5] S. Babaeinejadsarookolaee, B. Park, B. Lesieutre, and C. L. DeMarco, "Transmission congestion management via node-breaker topology control," *IEEE Systems Journal*, vol. 17, no. 3, pp. 3413–3424, 2023.

[6] Y. Jia, X. Wu, Z. Yang, and Y. Wang, "A two-stage approach for topology change-aware data-driven opf," *IEEE Transactions on Power Systems*, vol. 40, no. 2, pp. 1854–1867, 2025.

[7] M. Khanabadi, Y. Fu, and C. Liu, "Decentralized transmission line switching for congestion management of interconnected power systems," *IEEE Transactions on Power Systems*, vol. 33, no. 6, pp. 5902–5912, 2018.

[8] M. Heidarifar, P. Andrianesis, P. Ruiz, M. C. Caramanis, and I. C. Paschalidis, "An optimal transmission line switching and bus splitting heuristic incorporating ac and n-1 contingency constraints," *International Journal of Electrical Power & Energy Systems*, vol. 133, p. 107278, 2021.

[9] M. Kilthau, M. Asman, A. Karmann, G. Suriyamoorthy, J.-P. Beck, V. Regener, C. Derksen, N. Loose, M. Volkmann, S. Tripathi, F. Gehlhoff, K. Korotkiewicz, P. Steinbusch, V. Skwarek, M. Zdrallek, and A. Fay, "Integrating peer-to-peer energy trading and flexibility market with self-sovereign identity for decentralized energy dispatch and congestion management," *IEEE Access*, vol. 11, pp. 145 395–145 420, 2023.

[10] M. Attar, S. Repo, A. Mutanen, J. Rinta-Luoma, T. Väre, and K. Kukk, "Market integration and tso-dso coordination for viable market-based congestion management in power systems," *Applied Energy*, vol. 353, p. 122180, 2024.

[11] J. Wu, B. Zhang, Y. Jiang, P. Bie, and H. Li, "Chance-constrained stochastic congestion management of power systems considering uncertainty of wind power and demand side response," *International Journal of Electrical Power & Energy Systems*, vol. 107, pp. 703–714, 2019.

[12] H. Qiu, H. B. Gooi, H. Gao, V. Veerasamy, and Y. Li, "Cooperative power bidding of smart community grids with an aggregator-prosumer-based hierarchical framework," *IEEE Transactions on Smart Grid*, vol. 14, no. 4, pp. 3032–3045, 2023.

[13] "Pierre Pinson," [Online]. Available: https://pierrepinson.com/, accessed: Nov. 3, 2026.

[14] J. Viebahn, S. Kop, J. van Dijk, H. Budaya, M. Streefland, D. Barbieri, P. Champion, M. Jothy, V. Renault, and S. Tindemans, "Gridoptions tool: Real-world day-ahead congestion management using topological remedial actions," *CIGRE Science & Engineering*, no. C2.10528, 2024.

[15] A. Beiranvand and P. Cuffe, "Negative results on deploying distributed series reactance devices to improve power system robustness against cascading failures," *IEEE Transactions on Power Systems*, vol. 36, no. 6, pp. 5210–5221, 2021.

[16] MATLAB, *version 7.10.0 (R2019a)*. Natick, Massachusetts: The MathWorks Inc., 2019.

[17] R. D. Zimmerman, C. E. Murillo-Sánchez, and R. J. Thomas, "Matpower: Steady-state operations, planning, and analysis tools for power systems research and education," *IEEE Transactions on Power Systems*, vol. 26, no. 1, pp. 12–19, 2011.

[18] J. Lofberg, "YALMIP : a toolbox for modeling and optimization in MATLAB," in *2004 IEEE International Conference on Robotics and Automation (IEEE Cat. No.04CH37508)*, Sept 2004, pp. 284–289.

[19] A. Beiranvand, "Raw data and scripts from 'innovative integrated congestion management'," 2026. [Online]. Available: https://figshare.com/s/bb3f63922b6e69c345ff

[20] S. Babaeinejadsarookolaee, A. Birchfield, R. D. Christie, C. Coffrin, C. DeMarco, R. Diao, M. Ferris, S. Fliscounakis, S. Greene, R. Huang, C. Josz, R. Korab, B. Lesieutre, J. Maeght, T. W. K. Mak, D. K. Molzahn, T. J. Overbye, P. Panciatici, B. Park, J. Snodgrass, A. Tbaileh, P. V. Hentenryck, and R. Zimmerman, "The power grid library for benchmarking ac optimal power flow algorithms," 2021.

## Acknowledgment

ChatGPT (OpenAI GPT-5.5, Plus subscription) was used to assist with language refinement. All scientific, technical, and intellectual content, including the literature review, methodology, formulations, analyses, results, and conclusions, was conceived, developed, and critically reviewed by the author.

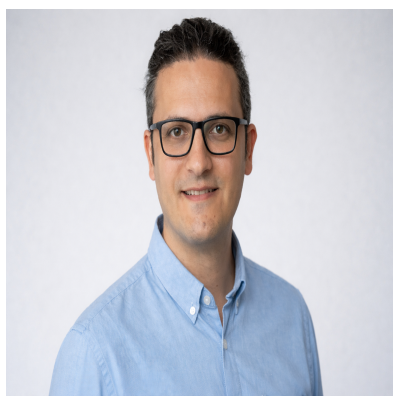

**Arash Beiranvand** (Senior Member, IEEE) is a Lecturer with the School of Electrical and Electronic Engineering, Technological University Dublin, Ireland. Previously, he was a Postdoctoral Research Fellow with the Energy Institute, University College Dublin, Ireland. He received the B.Sc. degree in Electrical Engineering (Power Systems) from Guilan University, Rasht, Iran, the M.Sc. degree in Electrical Engineering (Power Systems) from SBU, Tehran, Iran, and the Ph.D. degree in Electrical Engineering from University College Dublin, Ireland, in 2021. He is currently the Principal Investigator on several nationally funded research projects with a combined value exceeding €2.1 million. His research interests include power system operation, electricity markets, power system optimization, network topology optimization, renewable energy integration, and power system resilience.